\documentclass[12pt,twoside,final,psamsfonts]{amsart}
\usepackage[psamsfonts]{amssymb}
\usepackage{times,pictex,a4wide}
\theoremstyle{plain}         
\newtheorem*{teor*}{Theorem}
\newtheorem{teor}{Theorem}  
\newtheorem{prop}{Proposition}  
\newtheorem{lema}{Lemma}
\newtheorem{claim}{Claim}
\theoremstyle{remark}
\newtheorem*{remark}{Remark}
\newtheorem*{example}{Example}
\theoremstyle{definition}
\newtheorem*{defi}{Definition}
\newcommand{\C}{{\mathbb C}}
\newcommand{\N}{{\mathbb N}}
\newcommand{\D}{{\mathbb D}}
\newcommand{\R}{{\mathbb R}}
\renewcommand{\H}{{\mathcal H}}
\newcommand{\dbar}{\bar\partial}
\renewcommand{\Re}{\operatorname{Re}}
\renewcommand{\Im}{\operatorname{Im}}

\title{Multipliers and weighted $\dbar$ estimates}
\author{Joaquim Ortega-Cerd\`a}
\address{Dept.\ Matem\`atica Aplicada i An\`alisi,
Universitat  de Barcelona, Gran Via 585, 08071 Bar\-ce\-lo\-na, Spain}
\thanks{Supported by the DGICYT grant PB98-1242-C02-01, the CIRIT grant
2000-SGR00059 and the Departament d'Universitats i Recerca.
This is a corrected version of the paper published in Rev.\ Mat.\
Iberoamericana \textbf{18} (2002), 355--377. The typos and mistakes corrected here were
found during a formalization of the paper in Lean~4
(\texttt{github.com/joaquimortega/MR1949832}); the changes are kept to a
minimum}
\email{jortega@ub.edu}
\date{November 16 1999; corrected version November 2, 2026}

\begin{document} 

\begin{abstract} 
We study estimates for the solution of the equation $\dbar u=f$ in one
variable. The new ingredient is the use of holomorphic functions with precise
growth restrictions in the construction of explicit solutions to the equation.
\end{abstract}

\maketitle
\section{Introduction}  
In the present paper we will consider the equation $\dbar u =f$ in one
dimension. This equation plays a key role in the study of many problems in
complex analysis and, for this reason has been extensively studied. It is of
particular interest to have good estimates of the size of $u$ in terms of the
size of $f$ (see \cite{Berndtsson94} for a survey on the state of the art of
this problem). The purpose of this note is to show how a construction of
holomorphic functions with very precise growth restrictions can yield estimates
for the solutions to the $\dbar$-equation. With this tool we have been able to
obtain new proofs of some well-known results and some new estimates as well.

The most basic estimate is given by H\"ormander's theorem 
(see \cite[p.~92]{Hormander90a}):

\begin{teor*}[H\"ormander] 
Let $\phi$ be a subharmonic function defined in a
domain $\Omega\subseteq \C$ such that $\Delta \phi \ge \varepsilon$ for some
$\varepsilon>0$. Then there is  a solution $u$ to the equation $\dbar u=f$ such
that 
\[ \|ue^{-\phi}\|_2\lesssim \|fe^{-\phi}\|_2. 
\] 
\end{teor*}

\begin{remark} We write $f\lesssim g$ if there is a constant $K$ such that
$f\le K g$, and $f\simeq g$ if both $f\lesssim g$ and $g\lesssim f$.
Throughout, $\Delta=\partial^2/\partial x^2+\partial^2/\partial y^2=4\partial\dbar$
is the standard Laplacian, so that $\Delta\log|z|=2\pi\delta_0$.
\end{remark}

We will focus our attention on the case in which $\Omega$ is either the disk or
the whole plane. When $\Omega=\C$, M.~Christ has proved that the 
solution operator that solves the $\dbar$ equation with minimal weighted $L^2$
norm is also bounded  on weighted $L^p$ norms, where $1\le p\le \infty$ if we
assume some regularity on the weight (see \cite{Christ91}). His theorem is as
follows:

\begin{teor}[Christ]\label{Christ}
Let $\phi$ be a subharmonic function in $\C$ such that  there is a radius $r>0$ 
such that $\Delta\phi (D)\ge
1$ for any disk $D$ of radius $r$. Moreover we assume that $\Delta\phi$ is a
doubling measure. Then there is a  solution $u$ to the equation $\dbar u=f$
such that
\[
\|ue^{-\phi}\|_p\lesssim \|fe^{-\phi}\|_p,
\] 
for all $p\in[1,\infty]$. 
\end{teor}

As M.~Christ mentions, the doubling hypothesis on $\Delta \phi$ is not of an
essential  nature. It can be relaxed, but nevertheless one has to assume some
regularity on $\phi$ apart from the strict subharmonicity if one wants to
obtain $L^\infty$ estimates for instance. This is clearly seen in the following
example, due to Berndtsson:

\begin{example} 
Take $\phi(z)=\sum_{n\ge 3} \frac 1{n^2} \log|z-1/n|$. This is a  subharmonic
function in $\D$ that is bounded (above and below) in $1/2<|z|<1$ and moreover
$\phi(1/n)=-\infty$. Choose any smooth datum $f$ with support in a small disk
lying inside the annulus $1/2<|z|<1$ and such that $\int_{\D} f(z)\,dm(z)\ne 0$.

If there is a solution $u$ to the equation $\dbar u=f$ in $\D$ with 
 $\|ue^{-\phi} \|_\infty\lesssim \|fe^{-\phi}\|_\infty$, then 
$u(1/n)=0$ since the right-hand side is finite. In addition $u$ is
holomorphic outside the support of $f$. That means that $u$ is identically $0$
in a neighborhood of $\partial \D$. This cannot be so, 
because $0=\int_{\partial
\D} u\, dz= \int_{\D} \dbar u\, dm(z)\ne 0$.   
\end{example} 

There are more sophisticated examples due to Forn{\ae}ss and Sibony
\cite{ForSib89} that show that it is also impossible to have weighted $L^p$
estimates as in H\"ormander's theorem for any $p>2$ if we do not assume some
regularity on the weight.

In another direction, it is possible to extend H\"ormander's basic theorem to a 
larger class of weights including some non-subharmonic functions. This was done 
initially by Donnelly and Fefferman in \cite{DonFef83} and many others
afterwards (see \cite{BerCha99} and the references therein).  A variant of
their theorem  (in a particular case of a weight in the disk) is the following:

\begin{teor*}
Let $\phi$ be a subharmonic function in the unit disk $\D$ such that its
$(1-|z|^2)^2\Delta\phi>\varepsilon$ for some
$\varepsilon>0$.  Then there is  a solution $u$ to the equation $\dbar u=f$
with
\[
\int_{\D} \frac {|u(z)|^2 }{1-|z|^2}e^{-\phi}\, dm(z)
\lesssim
\int_{\D}  |f(z)|^2 e^{-\phi}(1-|z|^2)\, dm(z).
\]
\end{teor*}
For a simple proof of this case see \cite{BerOrt95}.

If we assume some regularity on the weight, we can extend this result
to $L^p$ norms. We require the Laplacian of the weight to be locally
doubling (see section~\ref{doubling} for the precise definition).
We will prove the following: 

\begin{teor}\label{meu}
Let $\phi$ be a subharmonic function in the unit disk $\D$ such that 
$\Delta\phi (D(z,r))> 1$ for some $r>0$ where $D(z,r)$ is
any hyperbolic disk with center $z\in \D$ and radius $r$. Moreover we assume
that $\Delta\phi$ is a locally
doubling measure with respect to hyperbolic
distance. Then there is  a solution $u$ to the equation $\dbar u=f$ with
\[
\int_{\D} \frac {|u(z)|^p}{1-|z|^2}e^{-p\phi}\, dm(z)
\lesssim
\int_{\D}  \frac{|f(z)(1-|z|^2)|^p}{1-|z|^2} e^{-p\phi}\, dm(z),
\]
for any $p\in[1,+\infty)$ and
\[
\sup |u| e^{-\phi} \lesssim \sup |f(\zeta)(1-|\zeta|)|e^{-\phi(\zeta)}.
\]
\end{teor}

\begin{remark} Observe that in the case $p\in[1,+\infty)$ we could have
rewritten the statement of the theorem if we absorb the factor $1/(1-|z|^2)$ in
the weight, i.e., if we replace $p\phi$ by $\widetilde\phi=p\phi+\log(1-|z|^2)$.
In this way it will look formally more similar to
H\"ormander's theorem, but we are allowing weights such that
$(1-|z|^2)^2\Delta\widetilde\phi>(-4+\varepsilon)$.
In particular, it includes functions $\widetilde\phi$
which are not subharmonic.

The weight $\phi/p$ satisfies the hypotheses of Theorem~\ref{meu} (with a
larger $r$). Applying the theorem to it we obtain, for each $p\in[1,+\infty)$,
a solution $u_p$ (which depends on $p$) with
\[
\int_{\D} \frac {|u_p(z)|^p}{1-|z|^2}e^{-\phi}\, dm(z)
\lesssim
\int_{\D}  \frac{|f(z)(1-|z|^2)|^p}{1-|z|^2} e^{-\phi}\, dm(z).
\]
\end{remark}

This is our main theorem, although the emphasis should be on the method of
proof rather than the new estimates. 
For instance, it is also possible to show
with the same type of proof that Theorem~\ref{Christ} holds when the measure
$\Delta \phi$ is locally doubling instead of doubling.

Our main tool (the multiplier) is a holomorphic function with very precise
growth restrictions. It is constructed in section~\ref{multiplier} and it may
exist under a less restrictive hypothesis, as in \cite{LyuMal99}. Our 
construction yields a more precise result that it is
needed when we
want to obtain estimates for the $\dbar$ equation.

With the same technique we can deal with some degenerate cases when
the weight $\phi$ is harmonic in large parts of the domain. In such a case one
has to impose extra conditions on the data of the equation, as in the following
theorem which may be of interest in the study of the so-called weighted
Paley-Wiener spaces.

\begin{defi} 
A positive Borel measure $\mu$ in $\C$ is a \emph{two-sided Carleson measure} 
whenever there is a
constant $C>0$ such that $\mu(D(x,r))\le C r$ for all disks of center
$x\in\R$ and any positive radius $r$.  
\end{defi}

\begin{teor}\label{degenerat}
Let $\phi$ be a subharmonic function in $\C$ such that the measure $\Delta
\phi$ is a locally 
doubling measure supported in the real line and  $\Delta \phi
(I(x,r))>1$ for some $r>0$ where $I(x,r)$ is the interval in $\R$ of center $x$
and radius $r$. Consider the equation $\dbar u=\mu$, where $\mu$ is a compactly
supported measure such that $e^{-\phi}d|\mu|$ is a two-sided  Carleson measure.
Then there is a solution $u$ with
\[
\limsup_{z\to\infty} |u(z)|e^{-\phi(z)}=0 \quad \text{ and } \quad
|u(x)|e^{-\phi(x)}\le C \left(1+ \int_{|z-x|<1}\frac {e^{-\phi(z)}\,d|\mu|(z)}{|x-z|}\right)
\]
for any $x\in \R$, where $C$ does not depend on the support of $\mu$.
\end{teor}

\begin{remark}
If $\phi$ is bounded below in the strip $|\Im z|<1$ (for instance if
$\phi=\alpha|\Im z|$), the factor $e^{-\phi(z)}$ in the integral can be
omitted. In general it cannot: the problem is invariant under the change
$\phi\mapsto\phi+\Re H$, $\mu\mapsto e^{H}\mu$, $u\mapsto ue^{H}$, with $H$
entire (which preserves the hypotheses and the left-hand side), while the integral
without this factor is not, and for the weight
$\phi(z)=|\Im z|+\Re z\,\Im z$ the estimate without the factor fails.
\end{remark}

The solution $u$ to the equation $f$ that we present is fairly explicit. It is
\emph{not} the canonical solution (i.e. the minimal $L^2$ weighted solution).
For instance in the case of Theorem~\ref{Christ} our
solution $u$ is given by an integral kernel
\begin{equation}\label{nucli}
u(z)=\int_\C e^{\phi(z)-\phi(\zeta)}k(z,\zeta) f(\zeta)\, dm(\zeta),
\end{equation}
which behaves differently from the canonical one. The kernel for the canonical
solution can sometimes be estimated. If the weight $\phi$ is of the form
$\phi(z)=b(x)$ and $0<c^{-1}<b''(x) <c$, then the kernel $k'$ of the canonical
solution has at most an exponential decay, i.e. there is a constant $A$ such
that $\limsup_{z\to\infty} |k'(z,0)|\exp(A|z|)=\infty$ 
(\cite[Proposition~1.18]{Christ91}). 
The kernel of our
solution has a much faster decay, namely 
\begin{prop}\label{churrimangui} 
Under the hypothesis of Theorem~\ref{Christ} there is a kernel
$k(z,\zeta)$ such that the function $u$ given by \eqref{nucli} is a solution to
the equation $\dbar u= f$ and for some $\varepsilon>0$,
\[
|k(z,\zeta)|\lesssim \frac{e^{-\varepsilon |z-\zeta|^2}}{|z-\zeta|}.
\]
\end{prop}
(A two-sided estimate cannot be expected: the kernel that we construct
vanishes when $z$ is a zero of one of the multipliers used to build it.)
However, there are some instances in which the canonical kernel has a faster
decay than our solution (when $\Delta \phi$ is very large). 

The structure of the paper is the following. In section~\ref{doubling} we will
prove some basic results on locally 
doubling measures which will be needed later. In
section~\ref{multiplier} we will construct our main technical tool, the
so-called multiplier. We will do so in the disk and in the whole plane. The
proof follows the same lines in both cases. 
Finally in
section~\ref{demostracions} we will show how we can use the multipliers to prove
Theorem~\ref{meu} and a new proof of 
Theorem~\ref{Christ} in which the doubling condition on $\Delta\phi$
is replaced by the locally
doubling condition. We will also sketch
how the same ideas can be used to prove Theorem~\ref{degenerat} and
Proposition~\ref{churrimangui}.

\section{Locally doubling measures}\label{doubling}
In this section we compile some basic facts we need on locally doubling 
measures. There are some intersections with the analysis of
\cite{Christ91}.  Recall that we always work in a domain $\Omega$ which is
either the plane or the disk. When the domain is $\C$ the natural distance is
Euclidean; in the case of $\D$ we will work with hyperbolic
distance. 

\begin{defi}
A measure $\mu$ in $\Omega$ is called a \emph{locally  doubling measure}
whenever there is a constant $C>1$  such that $\mu(B)\le C \mu(B')$, for all
balls $B\subset \Omega$ of radius smaller than $1$,  where $B$ is the ball
with the same center as $B'$ and twice its radius.  
\end{defi}

\begin{remark}
In the definition, we can replace the restriction that the radius of $B$ is
smaller than $1$ by any other constant $A$. The measures will be the same, but of
course the constant $C$ that appears will change. Indeed, if $1/2\le r<A$, the
ball $B(z,2r)$ is covered by a bounded number (depending on $A$) of balls of
radius $1/4$, and the center of each of them can be joined to $z$ by a chain
of at most $8A+1$ points at mutual distance less than $1/4$; applying the
definition along the chain we get $\mu(B(z,2r))\le C_A\,\mu(B(z,1/4))\le
C_A\,\mu(B(z,r))$. We use this freely below; in particular,
Lemmas~\ref{dis-mes} and~\ref{clau} hold with $1$ replaced by any constant
$A$, with constants depending on $A$.
\end{remark}

\begin{example}  
There are many locally doubling measures that are not doubling. They can grow
faster, for instance $d\mu(z)=e^{|z|}dm(z)$ is a locally doubling measure in
$\C$ equipped with Euclidean distance, while any doubling measure has at
most polynomial growth. Moreover they do not need to satisfy any strong 
symmetric
condition, for instance the measure $(\Im z)^3dm(z)$ for $\Im z>0$ and
$(\Im z)^2dm(z)$ for $\Im z<0$ is locally doubling and it is not doubling.
\end{example}

We start with an elementary lemma which is in fact an alternative 
description of locally doubling measures.

\begin{lema}\label{dis-mes}
Let $\mu$ be a locally
doubling measure in $\Omega$.  Then there is a $\gamma>0$
such that for any balls $B'\subset B$ of radius $r(B')$ and $r(B)<1$
respectively, we have
\[
\left(\frac{\mu(B)}{\mu(B')}\right)^{\gamma}\lesssim \frac {r(B)}{r(B')}
\lesssim \left(\frac{\mu(B)}{\mu(B')}\right)^{1/\gamma}.
\]
\end{lema}

\begin{proof}
The right inequality is essentially
Lemma~2.1 in \cite{Christ91} and the left inequality follows
directly from the definition. The converse is also true. If a measure
satisfies the inequalities with $B=2B'$ then it is locally doubling. 
\end{proof}

As a consequence of this lemma any locally doubling measure has no atoms. But it
is possible to prove more:
\begin{lema}\label{dos} 
Given any segment or circular arc $I\subset\Omega$ and any locally doubling measure
$\mu$ in $\Omega$, then $\mu(I)=0$.
\end{lema}
\begin{proof}
We consider the case of a segment; for an arc the same argument works in
polar coordinates centered at the center of the circle.
Assume that this is not the case. Then there is a subinterval $I'\subset I$ 
such that $\mu(I')>0$ and such that the square
of side length $|I'|$ that is bisected by $I'$ is inside $\Omega$ (see
figure~\ref{quadrat}). 
We define a  doubling 
measure $\nu$ in the interval $J$ which is the base
of the square that contains $I'$. The measure $\nu(A)$ of any set $A\subset J$ 
is defined as
$\nu(A)=\mu(R_A)$, where $R_A$ is the set in the square that projects orthogonally onto
$A$. Since $\mu$ is locally doubling, then $\nu$ is doubling for small
intervals: if $A$ is an interval, $R_{2A}$ is covered by $\simeq|J|/|A|$ balls
of radius $2|A|$ centered on the middle line of $R_A$ at mutual distance $|A|$,
and the concentric balls of radius $|A|/2$ are disjoint and contained in $R_A$.
Therefore $\nu$ has
no atoms. This implies that $\mu(I')=0$.
\end{proof}
\begin{figure}
\begin{center}
\font\thinlinefont=cmr5
\mbox{\beginpicture
\setcoordinatesystem units <1.04987cm,1.04987cm>
\unitlength=1.04987cm
\linethickness=1pt
\setplotsymbol ({\makebox(0,0)[l]{\tencirc\symbol{'160}}})
\setshadesymbol ({\thinlinefont .})
\setlinear
%
%
\linethickness= 0.500pt
\setplotsymbol ({\thinlinefont .})
\putrectangle corners at 13.636 21.882 and 16.161 19.357
%
%
\linethickness= 0.500pt
\setplotsymbol ({\thinlinefont .})
\setdashes < 0.1270cm>
\plot 14.882 21.880 14.897 19.363 /
\linethickness= 0.500pt
\setplotsymbol ({\thinlinefont .})
\setsolid
%
%
\plot 11.692 18.982 	11.800 18.885
	11.901 18.794
	11.995 18.711
	12.083 18.635
	12.165 18.565
	12.242 18.502
	12.314 18.444
	12.382 18.392
	12.506 18.304
	12.617 18.234
	12.719 18.180
	12.814 18.140
	12.916 18.108
	13.030 18.080
	13.155 18.056
	13.221 18.045
	13.290 18.035
	13.360 18.027
	13.431 18.019
	13.504 18.012
	13.578 18.005
	13.653 18.000
	13.728 17.995
	13.804 17.992
	13.881 17.989
	13.957 17.986
	14.033 17.985
	14.108 17.984
	14.183 17.984
	14.257 17.984
	14.330 17.985
	14.402 17.987
	14.471 17.990
	14.539 17.992
	14.605 17.996
	14.730 18.005
	14.844 18.015
	14.946 18.028
	15.067 18.045
	15.133 18.054
	15.202 18.065
	15.274 18.075
	15.349 18.087
	15.426 18.100
	15.506 18.113
	15.588 18.127
	15.671 18.143
	15.756 18.159
	15.841 18.177
	15.928 18.196
	16.015 18.216
	16.102 18.237
	16.189 18.260
	16.276 18.285
	16.363 18.311
	16.448 18.338
	16.533 18.367
	16.616 18.398
	16.697 18.431
	16.777 18.465
	16.854 18.502
	16.929 18.540
	17.001 18.581
	17.070 18.623
	17.135 18.668
	17.255 18.764
	17.359 18.870
	17.414 18.944
	17.463 19.026
	17.506 19.117
	17.544 19.213
	17.576 19.315
	17.602 19.421
	17.624 19.530
	17.641 19.640
	17.654 19.751
	17.662 19.861
	17.666 19.969
	17.667 20.075
	17.664 20.176
	17.659 20.271
	17.650 20.360
	17.638 20.441
	17.617 20.526
	17.584 20.615
	17.540 20.708
	17.487 20.804
	17.427 20.901
	17.362 21.001
	17.292 21.101
	17.221 21.202
	17.150 21.302
	17.080 21.402
	17.014 21.501
	16.953 21.598
	16.899 21.692
	16.853 21.783
	16.818 21.871
	16.796 21.954
	16.795 22.047
	16.820 22.153
	16.858 22.268
	16.901 22.388
	16.938 22.510
	16.958 22.631
	16.951 22.747
	16.908 22.854
	16.807 22.949
	16.744 22.977
	16.673 22.995
	16.598 23.006
	16.518 23.010
	16.437 23.011
	16.354 23.009
	16.272 23.009
	16.192 23.010
	16.116 23.017
	16.045 23.030
	15.925 23.084
	15.843 23.190
	15.830 23.259
	15.850 23.333
	15.909 23.420
	16.010 23.527
	/
%
%
\put{$I'$} [lB] at 15.058 20.665
%
%
\put{$J$} [lB] at 15.339 19.037
%
%
\put{$\Omega$} [lB] at 12.253 18.982
\linethickness=0pt
\putrectangle corners at 11.676 23.544 and 17.697 17.958
\endpicture}

\end{center}
\caption{}\label{quadrat}
\end{figure}
Let us introduce some notations. Let $\mu$ be a locally doubling measure in
$\Omega$ with $\mu(\Omega)=+\infty$.

\begin{defi}
For any $z\in\Omega$, denote by $\rho(z)$ the radius such that
$\mu(B(z,\rho(z))=1$. 
\end{defi}

This is always well defined since for any locally  doubling measure in
$\Omega$, the measure of any sphere is $0$ (with the same proof as in
Lemma~\ref{dos}). Thus the function $r\to
\mu(B(z,r))$ is continuous and strictly increasing. 

Since the measures that we consider (in Theorems \ref{Christ}, \ref{meu} and
\ref{degenerat}) are all measures such that
$\mu(B(z,r))\ge 1$ for some $r$ uniformly in $z$, then $\rho(z)$ has an upper
bound, but it can be very small.

The following claim 
is an immediate consequence of Lemma
\ref{dis-mes}.
\begin{claim}\label{propietats-mesura} 
Let $\mu$ be a locally doubling measure such that $\rho(z)$ 
has an upper bound. Then for any $K>0$ there is a $C_K$ such that 
$1/C_K<\rho(z)/\rho(w)<C_K$
whenever 
\[d(z,w)\le K \max(\rho(z),\rho(w)).\]
\end{claim}
Thus the radius of balls of measure one do not
change very abruptly. 
The following estimate is basic in our analysis:

\begin{lema}\label{clau}
If $\mu$ is a locally doubling measure in $\Omega$, then 
there is an $m\in\N$ such that for any  $\delta>0$ ,
\[
\sup_{w\in \Omega} \int_{\delta\rho(w)\le d(z,w)<1}
\left(\frac{\rho(z)}{d(z,w)}\right)^m \,d\mu(z)<C_\delta< +\infty.
\]
\end{lema}

\begin{proof}
We split the integral into two. In the first we integrate over the region
$\delta\rho(w)<\allowbreak d(z,w)<\rho(w)$. In this region $\rho(z)\simeq\rho(w)$,
therefore the integral is bounded by  
$C'_\delta\mu(B(w,\rho(w))$. In the second we integrate over the
region $\rho(w)\le d(z,w)\le 1$. 
We split this into rings of doubling size and we
may estimate it by
\[
\sum_{n=0}^k \int_{2^n<\frac{d(z,w)}{\rho(w)}<2^{n+1}}  \left(\frac{\rho(z)}
{2^n\rho(w)}\right)^m\, d\mu(z),
\]
where $k$ is such that $1<2^{k}\rho(w)\le 2$. 

Consider now the ball $B'$ of center $z$ and radius $\rho(z)$ and the ball $B$
of center $w$ and radius $C d(z,w)\simeq 2^n\rho(w)$. The constant $C$ is 
chosen
in such a way that $Cd(z,w)\ge \rho(z)+d(z,w)$. This is always possible,
since $\rho(z)$ and $\rho(w)$ are equivalent whenever $z$ is close to $w$.  
Therefore $B'\subset B$,
the radius of $B$ is smaller than $1$ and  we
may apply Lemma \ref{dis-mes}. We estimate $\rho(z)/(2^n\rho(w))$ by
$(C/\mu(B(w,2^n\rho(w)))^\gamma$, and the integral is bounded by a
constant times
\[
\sum_{n\ge 0}^k \frac{1}{(\mu(B(w,2^n\rho(w)))^{m\gamma-1}}=
\sum_{n\ge 0}^k \frac{\mu(B(w,\rho(w)))^{m\gamma -1}}
{(\mu(B(w,2^n\rho(w)))^{m\gamma-1}}.
\]

We apply Lemma \ref{dis-mes} to this expression and compare 
the quotient of measures by the quotient of radii (we consider
$1=\mu(B')=\mu(B(w,\rho(w))$ as the numerator) and obtain
\[
C\sum_{n=0}^k
\left(\frac{\rho(w)}{2^n\rho(w)}\right)^{\gamma(m\gamma-1)}<+\infty
\]
provided that we choose $m$ so that $m\gamma>1$.
\end{proof}

\section{The multipliers}\label{multiplier}
The main tool used to prove these results is the construction of the so-called
multipliers. These are holomorphic functions that have very precise growth
control. They have been used to solve some interpolation and sampling 
problems in several function spaces (see \cite{OrtSei98}, \cite{LyuSei94}) and
also the zero sets as in the Beurling-Malliavin theorem  (see also
\cite{Seip95b}). They all boil down to an approximation of subharmonic
functions by the logarithm of entire functions outside an exceptional set. 
The most general result of this
type is due to Lyubarski{\u\i} and Malinnikova, \cite{LyuMal99}, 
where they do not assume any
regularity condition on the Laplacian of the subharmonic function. However
we need a more precise description than theirs on 
the exceptional set in which the approximation need not hold.

The following theorem is a
result by Lyubarski{\u\i} and Sodin which will serve as a model
(see \cite{LyuSei94} for a proof).

\begin{teor*}[Lyubarski{\u\i}-Sodin]
Let $\phi$ be a subharmonic function in $\C$ such that its Laplacian 
satisfies $\Delta
\phi\simeq 1$. Then there exists an entire function $f$ with a uniformly 
separated 
zero set $Z(f)$ (i.e. $\inf_{a,b\in Z(f),a\ne b}d(a,b)>0$) such that 
\[
|f(z)|\simeq e^{\phi(z)},
\]
when $|z-a|\ge \varepsilon$ for all $a\in Z(f)$.
\end{teor*}

In the case of the disk the following theorem from Seip, \cite{Seip95b} is the
analogue to the multiplier lemma of Lyubarski{\u\i} and Sodin,

\begin{teor*}[Seip]
Let $\psi$ be a subharmonic function in $\D$ such that its Laplacian verifies
$(1-|z|^2)^2\Delta \psi\simeq 1$. Then there is a function $g\in\H(\D)$, with 
a  uniformly separated zero set $Z(g)$, and 
\[
|g(z)|\simeq e^{\psi(z)},
\]
when $\frac {|z-a|}{|1-\bar a z|}\ge \varepsilon$ for all $a\in Z(g)$.
\end{teor*}

We will need an analogous theorem for locally
doubling measures in $\C$ and in $\D$. 

\begin{teor}\label{multiplicador} 
Let $\psi$ be a subharmonic function in $\Omega$ such that its Laplacian
$\Delta \psi$ is a locally
doubling measure,  with the property $\Delta\psi (D(z,R))>1$
for all disks of some large radius $R>0$.  Then there is a holomorphic function
$h$ with zero set $Z(h)=\Lambda$ such that  
\[
\left(\frac{d(z,\Lambda)}{\rho(z)}\right)^M\lesssim
|h(z)|e^{-\psi(z)}\lesssim \frac{d(z,\Lambda)}{\rho(z)}
\]
for some fixed $M\in\N$, where $d(z,\Lambda)$ is the distance 
(in the appropriate metric) from $z$ to $\Lambda$.
\end{teor}
\begin{remark}
It will follow 
from the construction of $h$ that $d(z,\Lambda)\lesssim\rho(z)$,
thus  $|h|\simeq e^{\psi}$ 
outside a set $E_h$ composed  
of small disks around the zeros of $h$:
$E_h=\cup_{\lambda\in\Lambda}D(\lambda,\varepsilon\rho(\lambda))$.

With a slight refinement of the construction it is possible to prove that the
zero set $\Lambda$ can be chosen in such a way that $d(\lambda_i,\lambda_j)\ge
\varepsilon \max(\rho(\lambda_i),\rho(\lambda_j))$, for some $\varepsilon>0$
and $M$ can be chosen to be $1$, but we won't need this here.
\end{remark}
We will simultaneously prove Theorem~\ref{multiplicador} on the multipliers 
in the disk and in the plane, since we have to follow the same 
steps. To begin with,
we need a partition of the domain into rectangles that is 
well adapted to the measure
and the metric.  

\begin{lema}\label{partition}
Assume that $\mu$ is a locally doubling measure in $\Omega$ with
$\mu(D(z,R))>1$.
Given any $N\in \N$ there is  a partition of $\Omega$ into rectangles
$\{R_i\}_{i\in I}$ such  that  $\mu(R_i)=N$ and if we denote by $L_i$
the length of the longer side of $R_i$ and $l_i$ the 
length of the smaller side,
then  $sup_{i\in I} L_i/l_i=L(N)<+\infty$. 
\end{lema}

\begin{remark} When $\Omega$ is a disk, one has to understand that
by ``rectangles''
we mean preimages, under the Cayley transform $w=i(1+z)/(1-z)$ of $\D$ onto
the upper half-plane, of rectangles with sides parallel to the axes, i.e.,
regions bounded by two horocycles and two geodesics ending at $1$; the side
lengths are those of the rectangles in the upper half-plane. This lemma is basically
the partition theorem from \cite{Yulmukhametov85}, but we include a proof, since
the doubling assumption (which is not needed) makes it particularly easy. 
\end{remark}

\begin{proof}
We start by assuming that $N=1$, the general case follows if we use the same
construction with the measure $\sigma=\mu/N$ instead of the measure $\mu$. We
will first find a partition into rectangles $\{\widetilde R_i\}_{i\in I}$ in
such a way that $\mu(\widetilde R_i)\in \N$, $1\le\mu(\widetilde R_i)\le C$ 
and with the ratio of side-lengths bounded
and $C$ is a constant that depends only on the doubling constant (and on
$R$). Later on, we will refine this
partition in order to obtain rectangles of mass one.

Recall that there is some $R>0$ such that $\mu(D(z,R))>1$ for all $z\in\Omega$.
Let us partition the plane (see next paragraph for the disk)
into parallel strips of width $2R$. Then, we
slice each strip in rectangles of mass a natural number (the sides of the
rectangle have no mass because of Lemma~\ref{dos}).
The length of any piece will be between $2R$ and $4R$.
Since any square of size $2R\times 2R$ contains a disk of radius $R$, it has
mass at
least $1$, and it is possible to slice the strip in such a way that the
resulting rectangles have the ratio between the sides bounded by $2$.
We have no upper bound of the mass of these rectangles;
we only know that it is a natural number.

In the case of the domain being the disk, we work in the upper half-plane,
the image of $\D$ by the Cayley transform, where the hyperbolic metric is
$|dw|/\Im w$. We replace the strips by the horizontal strips
$\Lambda^j\le \Im w<\Lambda^{j+1}$, $j\in\mathbb Z$, where $\Lambda>1$ is a large
constant depending on $R$, and we slice each of them in rectangles of
integer mass, with horizontal side between $h_j$ and $2h_j$, where
$h_j=(\Lambda-1)\Lambda^j$ is the height of the strip. If $\Lambda$ is large
enough, any square of side $h_j$ in the strip contains a hyperbolic disk of
radius $R$, so this can be done as in the plane.
The resulting rectangles have the ratio of lengths of
sides bounded by $2$ and bounded hyperbolic diameter (depending on
$\Lambda$), and again there is no control on the upper bound of the mass. In
what follows, ``disks'' are hyperbolic disks of bounded radius, to which the
doubling property applies by the Remark in Section~\ref{doubling}.

From now on the procedure in the disk and in the plane will be the same. 
We will divide each rectangle in two. All the resulting rectangles will still
have
integer mass and the ratio of the sides will always remain
bounded by $3$. We will
bisect each rectangle until the mass is smaller than the
constant $C$.

The bisection is done as follows: denote the original rectangle by
$R=[a,a+w]\times [b,b+l]$ where $l\le w\le 3l$. Consider a smaller auxiliary
 rectangle
$R'\subset R$ of the form $R'=[a+w/2-h/2,a+w/2+h/2]\times[b,b+l]$ with an 
auxiliary $h\le w$ so
that $\mu(R')=1$. The length $h$ of $R'$ cannot be very large. 
If $h>w/3$, then
$R\subset 3R'$ (the rectangle with the same center as $R'$ and sides three
times as long), and $R'$ has ratio of side-lengths bounded by $3$. Thus $R'$
contains a disk and $R$ is contained in a concentric disk of comparable
radius, so that $\mu(R)\le C\mu(R')= C$, where $C$ depends only on the
doubling constant, and we do not need to bisect $R$. So we may assume that
$h\le w/3$. Since $\mu(R')=1$, there must be a $0<t<h$ such that
$\mu(R_1)\in\N$, where $R_1=[a,a+w/2-h/2+t]\times[b,b+l]$. We
denote by $R_2=R\setminus R_1$. It also verifies $\mu(R_2)\in \N$. Both
masses are positive, because a locally doubling measure with
$\mu(\Omega)=+\infty$ gives positive mass to every disk, hence to every
rectangle. Finally it is
easily checked that $l/3\le w/3\le w/2-h/2+t\le w\le 3l$, therefore the
quotients of the sides-length of $R_1$ are bounded by 3 and similarly with
$R_2$.

Thus far, the rectangles are not very deformed and all have a mass between
$1$ and $C$. In order to obtain rectangles of mass $1$, we 
divide each of them into
rectangles of mass one by cutting along the direction of the longest side. 
The local
doubling condition ensures that all of them
will be essentially of the same proportion. Indeed, let $S$ be a slab of
width $s\le l$ of a rectangle of sides $w\ge l$, cut along the longest side.
Cover $S$ by $\simeq l/s$ disks of radius $s$ centered on its middle line, and
compare each of them, by the doubling property, with the disk of radius $s/2$
whose center is at distance $2s$ along the longest side; we obtain
$\mu(S)\le C'\mu(S')$, where $S'$ is the adjacent slab of width $2s$ and $C'$
depends only on the doubling constant. Iterating, the slabs of width $3^js$
starting at $S$ (in either direction) have mass at least $(1+1/C')^j\mu(S)$. Since
the rectangle has mass at most $C$, a slab of mass one has width at least
$c\,w$, where $c>0$ depends only on $C$ and $C'$.
At most we are dividing
each rectangle in $C$ pieces, therefore the resulting rectangles have 
a bounded ratio of
side-lengths as desired.
\end{proof}
The family of rectangles that we have just constructed looks very much like
squares, since the excentricity is bounded, but moreover the size of the
rectangles changes very slowly, along with $\rho(z)$:

\begin{claim}\label{propietats} The family of rectangles $\{R_i\}$ constructed
in Lemma~\ref{partition} has the following two properties:
\begin{itemize}
\item The ratio between the diameter of $R$ and $\rho(z)$ for any $z\in R$ is
bounded above and below by two constants independent of $R$ and $z\in R$.
\item For any $K>0$ there is a constant $C_K>0$  such that whenever $K R_i\cap
K R_j\neq \emptyset$ the ratio between the diameters of $R_i$ and $R_j$ is
bounded by $C_K$. 
\end{itemize} 
\end{claim}
\begin{proof}
The first assertion follows since $R$ has bounded excentricity and constant
mass. The second one is an immediate consequence of 
Claim~\ref{propietats-mesura}.
\end{proof}
In order to construct the multiplier, we first select  its zeros.
We take a very large $N=mk$ (the same $m$ as given by Lemma~\ref{clau}, which
we may take $\ge 3$, and
an even $k\ge m$, as required in Lemma~\ref{moments}).  We partition
$\Omega$ in rectangles $\{R_i\}_{i\in I}$ of mass $N$ for the measure
$\Delta\psi/(2\pi)$ as in
Lemma~\ref{partition} (since $\Delta\log|z|=2\pi\delta_0$, a zero of $h$
accounts for a mass $2\pi$ of $\Delta\psi$). For any $i\in I$, we will choose $N$ points
$\{\lambda^i_1,\ldots,\lambda^i_N\}$ which lie near $R_i$ and such that the
moments of order $0,1,2,\ldots m-1$  of the measure $\Delta\psi/(2\pi)$ restricted to
$R_i$ coincide with the corresponding moments of the measure $\sum_{j=1}^N
\delta_{\lambda^i_j}$. The following lemma addresses this point.

\begin{lema}\label{moments}
Let $D$ be a closed disk in $\C$.
Given any $m\in \N$, any $k\in\N$ with $k\ge m$ and any $C>1$,
for any measure $\mu$  supported in $D$ of  total mass $N=mk$, there are two  sets of $N$
points  $\Lambda=\{\lambda_1,\ldots,\lambda_N\}$ inside $2mD$
and $\kappa=\{\kappa_1,\ldots,\kappa_N\}$ inside $(C+4m)D\setminus CD$
satisfying
\[
\int_D z^j\, d\mu(z)=
\lambda_1^j+\cdots +\lambda_N^j=\kappa_1^j+\cdots+\kappa_N^j ,\qquad j=0,\ldots,m-1.
\]
For a rectangle $R$ we apply it to $D=D_R$, the smallest closed disk
containing $R$, and we write $\Lambda(R)$ and $\kappa(R)$.
\end{lema}

\begin{proof}
We want that
\[
\int_D p(z)\, d\mu(z)=\sum_{i=1}^N p(\lambda_i),
\]
for all polynomials of degree less than or equal to $m-1$.
Equal-weight (Chebyshev) quadrature formulas with nodes inside $D$ (see
\cite{Gautschi76} or \cite{Korevaar94}) need not exist for an arbitrary
measure (for instance, for $\mu=N(w\delta_{-1}+(1-w)\delta_1)$ with $w$
irrational and $m=3$ there are no such nodes in the square $[-1,1]^2$), so we
allow the nodes to lie in $2mD$. Let $c$ and $r$ be the center and the radius of
$D$, and let $s_j=\frac 1k\int_D (z-c)^j\,d\mu(z)$, so that $|s_j|\le m r^j$,
$j=1,\ldots,m$. Define $e_0=1$ and, recursively,
$je_j=\sum_{i=1}^{j}(-1)^{i-1}e_{j-i}s_i$; by induction $|e_j|\le (mr)^j$.
Let $z_1,\ldots,z_m$ be the roots of $z^m-e_1z^{m-1}+\cdots+(-1)^me_m$. By
Newton's identities $z_1^j+\cdots+z_m^j=s_j$ for $j=1,\ldots,m$, and
$|z_i|\le 2mr$, since $|z|>2mr$ implies
$|z|^m>\sum_{j=1}^m (mr)^j|z|^{m-j}$.
We take as $\Lambda$ the points $c+z_1,\ldots,c+z_m$, each of them repeated
$k$ times. Then $\sum_{\lambda\in\Lambda}(\lambda-c)^j=ks_j=\int_D (z-c)^j\,d\mu(z)$
for $j=1,\ldots,m-1$, and also for $j=0$ because $\mu(D)=mk$.
These are the points that will be used in the construction
of the multiplier; they will be in fact the zeros of it. Note that all the
points  $\lambda_j$ appear with a multiplicity $k$ since there are $N=km$
points with
equal weights. For later use, it is convenient to have an alternative set of
zeros $\kappa_1,\ldots,\kappa_N$  at our disposal which are separated from the original
ones (outside $CD$) and still have the same moments.
This is easily done. It can be checked immediately that
$kp(\lambda)=\sum_{l=0}^{k-1}p(\lambda+ \tau e^{l2\pi i/k})$, for any
$\tau\in
\C$ and any polynomial of degree less than $k$, in particular of degree
$m-1$. Thus, we could take as an alternative
set  $\kappa_{j,l}=c+z_j+\tau e^{l2\pi i/k}$, $j=1,\ldots,m$, $l=0,\ldots,k-1$,
where $\tau=(C+2m)r$, so that all $\kappa_{j,l}$ are outside $CD$ and inside $(C+4m)D$.
\end{proof}

Now we take a holomorphic function $h$ that vanishes at all the points 
$\{\lambda_j^i\}_{i\in I,j=1,\ldots,N}$, where
$\{\lambda_j^i\}_{j=1,\ldots,N}=\Lambda(R_i)$ is given by Lemma~\ref{moments}
applied to the restriction of $\Delta\psi/(2\pi)$ to $R_i$ (the rectangles defined
in Lemma~\ref{partition}; see below for a modification in the disk).
This function is defined up to a
factor of
the form $e^g$, with $g\in \H(\Omega)$. We choose this $g$ in such a way that
\[
\log|h|= \psi - \frac 1{2\pi} \int_{\C}
\log|z-\zeta| \Bigl(\Delta\psi-2\pi\sum \delta_{\lambda_j^i}\Bigr),
\]
in the case of $\Omega=\C$ and
\[
\log|h|=\psi -\frac 1{2\pi}
\int_{\D} \log\left|\frac{z-\zeta}{1-\bar\zeta z}\right|
\Bigl(\Delta\psi-2\pi\sum \delta_{\lambda_j^i}\Bigr)
\]
in the case of $\Omega=\D$. In both cases $\psi$ is the subharmonic function in
the statement of Theorem~\ref{multiplicador}.
Thus the problem has been reduced to show that
\begin{equation}\label{diferenciaC}
M\log{\frac{d(z,\Lambda)}{\rho(z)}}-C\le
-\frac 1{2\pi}\int_{\C} \log|z-\zeta| \Bigl(\Delta\psi-2\pi\sum \delta_{\lambda_j^i}\Bigr)\le
\log{\frac{d(z,\Lambda)}{\rho(z)}}+C,
\end{equation}
in the case of $\Omega=\C$, and  when $\Omega=\D$, we have to obtain
\begin{equation}\label{diferenciaD}
M\log{\frac{d(z,\Lambda)}{\rho(z)}}-C\le
-\frac 1{2\pi}\int_{\D} \log\left|\frac{z-\zeta}{1-\bar\zeta z}\right|
\Bigl(\Delta\psi-2\pi\sum \delta_{\lambda_j^i}\Bigr)\le
\log{\frac{d(z,\Lambda)}{\rho(z)}}+C.
\end{equation}

The integral in \eqref{diferenciaC} is split as
\[
\sum_{i\in I}\int_{\C} \log |z-\zeta|
\left(\chi_{R_i}(\zeta)\Delta\psi(\zeta)-2\pi\sum_{j=1}^N
\delta_{\lambda_j^i}(\zeta)\right).
\]

In any of these integrals we can subtract any polynomial of degree $m-1$ to the
logarithm since the moments up to order $m-1$ of
$\chi_{R_i}(\zeta)\Delta\psi(\zeta)$ and $2\pi\sum_{j=1}^N  \delta_{\lambda_j^i}
(\zeta)$
are the same. Since the points $\lambda^i_j$ lie near $R_i$, but not
necessarily inside it, we say that two rectangles $R_i$, $R_j$ are
\emph{neighbors} if $d(R_i,R_j)<K\max(\operatorname{diam}R_i,\operatorname{diam}R_j)$,
where $K\ge 2m+2$ is a large constant. Then the points attached to a rectangle
are at distance at least $\operatorname{diam}R_i$ from any rectangle $R_i$ which is
not one of its neighbors, and by Claim~\ref{propietats} each rectangle has a
bounded number of neighbors, all of diameter comparable to its own.
For any $R_i$ far from $z$ (we exclude the rectangle  $R_j$ to
which $z$ belongs and its neighbors) we take a polynomial
$p$  of degree
$m-1$, which is the Taylor expansion of $\log|z-\zeta|$ about some point 
$\lambda_0^i\in R_i$.

The difference between $|\log|z-\zeta|-p(\zeta)|$ is bounded by $\frac
C{|z-w|^{m}}|\zeta-\lambda_0^i|^m$, where $w$ is some point in the segment
joining $\lambda_0^i$ and $\zeta$. Since $z$
does not belong to $R_i$ or any of its neighbors and $\zeta\in R_i$ (or
$\zeta$ is one of the points $\lambda^i_j$),
it follows that $|z-w|\simeq |z-\zeta|$ and  $|\zeta-\lambda_0^i|\lesssim
\rho(\zeta)$ by Claim~\ref{propietats}. Thus the integral is bounded by a constant times
\[
\int_{R_i}\frac{\rho(\zeta)^{m}}{|z-\zeta|^{m}}\Delta\psi(\zeta)+
2\pi N\frac{\rho(\lambda_0^i)^{m}}{|z-\lambda_0^i|^{m}}.
\]

Both the integral and the sum are comparable since $\rho(\zeta)\simeq
\rho(\lambda_0^i)$, $|z-\zeta|\simeq |z-\lambda_0^i|$ and the
$\Delta\psi$-mass of the rectangle is $2\pi N$.
This estimate is true for all $R_i$ except the one that contains $z$
and its neighbors. There is a $\delta>0$ (uniform in $z$ because of
Claim~\ref{propietats-mesura}), such that the sum over 
all such rectangles is bounded by
\[
\int_{|z-\zeta|\ge \delta\rho(z)} 
\frac{\rho(\zeta)^m}{|z-\zeta|^{m}}\Delta\psi(\zeta).
\] 
Call $M'=2\sup_z \rho(z)$.
If we integrate in the region $\delta\rho(z)\le|z-\zeta|\le M'$ we may apply
Lemma~\ref{clau} (with $1$ replaced by $M'$, see the Remark in
Section~\ref{doubling}). If we integrate in the region $|z-\zeta|\ge M'$, we may
estimate the integral by
\[
\int_{|z-\zeta|\ge M'}
K\frac{\rho(\zeta)^2}{|z-\zeta|^{3}}\Delta\psi(\zeta).
\]
We use that $\rho(\zeta)^2\simeq \int_{|\zeta-w|\le\rho(\zeta)}dm(w)$
and then
Fubini's theorem to obtain
\[
\int_{|z-w|\ge M'/2}
K\frac{1}{|z-w|^{3}}\, dm(w)<+\infty.
\]

There are at most a finite number of neighboring rectangles (uniformly
in $z\in \C$) to the rectangle that contains $z$ because all of them have size
comparable to $\rho(z)$. In each of them, since $\Delta\psi(R_i)=2\pi N$,
the contribution of $R_i$ and its points $\lambda_j=\lambda^i_j$ to the
middle term of \eqref{diferenciaC} is
\[
-\frac 1{2\pi}\int_{R_i} \log \frac{|z-\zeta|}{\rho(z)}\Delta\psi(\zeta) +
\sum_{j=1}^N\log\frac{|\lambda_j-z|}{\rho(z)}.
\]
The integral is bounded whenever $\Delta\psi$ is locally doubling. This is
Lemma~2.3 of \cite{Christ91} which is in turn a  direct consequence of
Lemma~\ref{dis-mes}. The sum accounts for the terms
$\left(\frac{d(z,\Lambda)}{\rho(z)}\right)^M$ and $\frac{d(z,\Lambda)}{\rho(z)}$
in the statement of the theorem: each of its boundedly many terms lies between
$\log\frac{d(z,\Lambda)}{\rho(z)}$ and a constant, and if
$d(z,\Lambda)<\rho(z)/C$ the point of $\Lambda$ nearest to $z$ is attached to
one of these rectangles, so that one of the terms is
$\log\frac{d(z,\Lambda)}{\rho(z)}$. Here $M$ is the maximal number of points
attached to the neighbors of a rectangle.

In the disk the points $\lambda^i_j$ must be chosen with some care. When
$\Delta\psi$ is sparse, the rectangles $R_i$ are hyperbolically large, and the
points given by Lemma~\ref{moments} (applied to the Euclidean disk $D_{R_i}$)
may fall outside $\D$. We fix a small $\rho_0>0$, depending on $m$ and on the
constants in Lemma~\ref{partition} and Claim~\ref{propietats}. If $\rho\le\rho_0$
on $R_i$, then $D_{R_i}$ has radius $\lesssim\rho_0(1-|\zeta|^2)$ for
$\zeta\in R_i$, so the points of Lemma~\ref{moments} lie in $\D$, at
hyperbolic distance $\lesssim\rho(\zeta)$ from $R_i$, and we use them. Otherwise we
take as $\Lambda(R_i)$ the barycenter
$b_i=N^{-1}\int_{R_i}\zeta\,\Delta\psi(\zeta)/(2\pi)$ repeated $N$ times, and as
$\kappa(R_i)$ (when it is needed in Proposition~\ref{tecnic}) the points
$b_i\pm v_i$, each repeated $N/2$ times, where $v_i$ is orthogonal to $b_i$ and
the hyperbolic distance from $b_i$ to $b_i\pm v_i$ is a large multiple of
$\rho(b_i)$. The barycenter lies in the hyperbolic convex hull of $R_i$, hence
in $\D$, but for these rectangles only the moments of order $0$ and $1$ are
matched. This is enough: a rectangle with $\rho>\rho_0$ which is not a neighbor
of the rectangle that contains $z$ lies at hyperbolic distance at least
$\delta_0$ from $z$, for some $\delta_0\in(0,1]$ depending on $\rho_0$, and in that region the
estimate below only uses the moments of order $0$ and $1$ (the Taylor formula
being applied along segments contained in a hyperbolic disk that contains $R_i$
and its points).

We will to estimate now the integral in \eqref{diferenciaD}, which can be
expressed as
\[
\sum_{i\in I}\int_{\D} \log \frac{|z-\zeta|}{|1-\bar\zeta z|}
\left(\chi_{R_i}(\zeta)\Delta\psi(\zeta)-2\pi\sum_{j=1}^N
\delta_{\lambda_j^i}(\zeta)\right).
\]

As before we can subtract a Taylor polynomial of degree $m-1$
(of degree $1$ for the large rectangles)
at a point $\lambda_0^i\in R_i$. Now, since
\[
\left|\nabla_\zeta^m \log \frac{|z-\zeta|}{|1-\bar\zeta z|}\right|\lesssim
\frac{1-|z|^2}{|1-\bar\zeta z||z-\zeta|^m},
\]
the integral is bounded by
\begin{equation}\label{disc}
C\int_{\zeta\notin\delta D(z,\rho(z))} 
\frac{(1-|z|^2)(1-|\zeta|^2)^m\rho(\zeta)^m}
{|1-\bar\zeta z||z-\zeta|^m}\Delta\psi(\zeta)+\sum
\log\frac{|z-\lambda_i|}{\rho(z)|1-\bar\lambda_i z|},
\end{equation}
where the sum is over all $\lambda_i$ that are attached to the rectangle $R_i$ which
contains $z$ and its neighbors.

We
split the integral in two pieces. In the first we integrate over the domain
$\Omega_1=\{\zeta\in \D;\ d(z,\zeta)<\delta_0, \zeta\notin\delta
D(z,\rho(z))\}$ (among the rectangles which are not neighbors of the one
containing $z$, only those with $\rho\le\rho_0$ meet $\Omega_1$),  and we use Lemma~\ref{clau} to obtain
\[
\int_{\Omega_1}\frac{\rho(\zeta)^m}{d(z,\zeta)^m}\Delta\psi(\zeta)<\infty.
\]

In $\Omega_2$ we have that $d(z,\zeta)\ge\delta_0$, and \eqref{disc} (or its
analogue with $m=2$ for the large rectangles)
is bounded by a constant times
\[
\int_{\Omega_2}\frac{(1-|z|^2) (1-|\zeta|^2)^2\rho(\zeta)^2}{|1-\bar\zeta z|^3}
\Delta\psi(\zeta).
\]

We may think of $(1-|\zeta|^2)^2\rho(\zeta)^2$ as
$\int_{d(w,\zeta)<\rho(\zeta)}dm(w)$ and apply Fubini's theorem to obtain 
\[
\int_{\Omega_2}\frac{(1-|z|^2) (1-|\zeta|^2)^2\rho(\zeta)^2}{|1-\bar\zeta z|^3}
\Delta\psi(\zeta)\le \int_{\D} \frac{(1-|z|^2)}{|1-z\bar w|^3}\, dm(w)<+\infty.
\qed
\]
Theorem~\ref{multiplicador} is not yet what we need for the estimates to the
$\bar\partial$-equation because the exceptional set of the multiplier introduces
a technical difficulty. We need to approximate the weight by a holomorphic
function \emph{everywhere}. This obstruction 
can be avoided using several multipliers
simultaneously as
described in the next proposition:
\begin{prop}\label{tecnic} Given $\psi$ as in the statement of
Theorem~\ref{multiplicador} there is a
collection of multipliers $h_1,\ldots,h_n$ satisfying the conclusion of
Theorem~\ref{multiplicador}. Moreover their exceptional sets 
(see the remark after
Theorem~\ref{multiplicador}) are disjoint; i.e.
$E_{h_1}\cap\cdots\cap E_{h_n}=\emptyset$. 
\end{prop}
\begin{proof}
Take the partition of $\Omega$ in rectangles given by Lemma~\ref{partition}. We
distribute the rectangles in a finite number of families of rectangles
$\Omega=\cup_{l=1}^n(\cup_{i\in I_l}R^l_i)$ with the property that any two
rectangles of the same family $R_i^l$, $R_j^l$ are very far apart (i.e.
$M R_i^l\cap M R_j^l=\emptyset$, for some large constant $M$). This is possible
with the Besicovitch covering lemma. Now for each family $\{R^l_i\}_{i\in I_l}$
we can construct a multiplier $h_l$ in such a way that it has no zeros in any of
the rectangles of the family $R^l_i$ nor in their neighbors. 
The way to proceed to construct $h_l$ is the following: For any rectangle
$R$ 
that is neither from the family $\{R^l_i\}_{i\in I_l}$ nor one of their
neighbors we take the set of points $\Lambda(R)$ given by Lemma~\ref{moments}.
For the rectangles $R$ from the family or their neighbors we use the
alternative set of points $\kappa(R)$ also defined in Lemma~\ref{moments}.
We build as before a multiplier $h_l$ with zeros at the selected points. 
It has the right growth and the
additional property that it has no zeros in the rectangles from the family 
$\{R^l_i\}_{i\in I_l}$ and
their neighbors. This is clear because we can choose a constant $C$ in
Lemma~\ref{moments} in such a way that the points $\kappa(R)$ are neither in $R$
nor in its neighbors. 
Moreover they are not so far apart from $R$ that they reach another
rectangle from the family 
(this can be prevented by choosing a very large $M$ in
the splitting of the rectangles into families). Thus the exceptional set for
$h_l$ does not include any rectangle from the family $\{R^l_i\}_{i\in I_l}$.
\end{proof}
\section{The $\bar\partial$-estimates}\label{demostracions}
This section contains three parts. In the first one, we will see how the
weights that we consider can be regularized without losing generality. In the
second subsection we prove the $L^p$ weighted $\bar\partial$-estimates in the
plane and the disk. Finally in the last part we indicate how
Theorem~\ref{degenerat} can be proved.

\subsection{The regularization of $\phi$}
In the hypothesis of the theorem we assume that for 
some large radius $r>0$, $\Delta \phi(D(z,r))>1$  at any point $z\in\Omega$. 
This is a condition that 
ensures that $\phi$ is ``strictly subharmonic''. It will be more convenient for
us to assume that $\Delta\phi >\varepsilon dm(z)$. This means that the measure
is  more regular since there are no ``holes'' with zero measure. The following 
proposition allows us to do so:

\begin{lema}\label{regularitzacio}
If the measure $\Delta\phi$ is a locally 
doubling measure in $\Omega$  and $\Delta
\phi(D(z,r))>1$ for  some large radius $r>0$ and any point $z\in\Omega$ then
there is a subharmonic  weight $\psi$ equivalent to the original, i.e.
$\sup_\Omega |\phi-\psi|<+\infty$,  such that $\Delta \psi$ is a locally 
doubling
measure and moreover  $\Delta\psi >\varepsilon dm(z)$ for some $\varepsilon>0$
(when $\Omega=\D$, $(1-|z|^2)^2\Delta\psi>\varepsilon$).
\end{lema}

\begin{proof}
We will split $\Delta \phi$ in two measures $\mu_1+\mu_2$.   To
describe the measure $\mu_1$,  let us tile the plane into squares
$Q_j$ of diameter $R>0$ (in the case of the disk, the pieces of a partition
of $\D$ as in Lemma~\ref{partition} for the hyperbolic area measure, with a
large mass, which have hyperbolic diameter comparable to a large constant $R$) in such a way
that $\Delta\phi(Q_j)>2$ for all  $Q_j$. This is feasible because of the
hypothesis on the measure. The measure  $\mu_1$ is defined as
$\mu_1|_{Q_j}=\frac 1{\Delta\phi(Q_j)}\Delta\phi$.  The measure $\mu_2$ is
$\Delta\phi-\mu_1$. 
It follows from the definition that  $\frac 12 \Delta\phi \le \mu_2\le
\Delta\phi$, therefore $\mu_2$ is a locally
doubling  measure. It is also true that $\mu_1$ is locally doubling
 because $\Delta\phi(Q_j)$ does not change abruptly in neighboring
squares
 and moreover $\mu_1(Q_j)=1$.  

We will regularize the measure $\mu_1$ by taking the convolution  (the
invariant convolution when $\Omega$ is a disk) of it  with the normalized
characteristic
function of a very large disk: $\widetilde \mu_1= \mu_1 \star
\frac{\chi_{D(0,2R)}}{|D(0,2R)|}$. The measure  $\widetilde\mu_1$ in the plane 
satisfies
$\varepsilon dm(z) <\widetilde \mu_1<K dm(z)$ (when $\Omega=\D$, it satisfies
$\varepsilon<(1-|z|^2)^2\widetilde\mu_1< K$). 

Instead of solving the Poisson equation for $\mu_1$ and $\widetilde\mu_1$
separately with global potentials (the corrections that are needed to make
such potentials converge do not commute with the convolution), we use a local
potential. For $\zeta\in\Omega$ let $\sigma_\zeta$ be the normalized area
measure on the disk $D(\zeta,2R)$ (the normalized hyperbolic area when
$\Omega=\D$), so that $\widetilde\mu_1=\int_\Omega\sigma_\zeta\,d\mu_1(\zeta)$,
and let
\[
G(z,\zeta)=\frac 1{2\pi}\left(\int_\Omega \log[z,w]\,d\sigma_\zeta(w)-\log[z,\zeta]\right),
\]
where $[z,w]=|z-w|$ if $\Omega=\C$ and $[z,w]=|z-w|/|1-\bar wz|$ if
$\Omega=\D$. Since $\Delta_z\log[z,w]=2\pi\delta_w$, we have
$\Delta_zG(\cdot,\zeta)=\sigma_\zeta-\delta_\zeta$. Since $w\mapsto\log[z,w]$ is
harmonic in $\Omega\setminus\{z\}$ and $\sigma_\zeta$ is radial with respect
to $\zeta$ (after a M\"obius transformation that maps $\zeta$ to $0$, in the
case of the disk), the mean value property gives $G(z,\zeta)=0$ when
$d(z,\zeta)\ge 2R$, while $|G(z,\zeta)|\lesssim 1+\log\frac{2R}{d(z,\zeta)}$
when $d(z,\zeta)<2R$.
 
We take as $\psi=\phi +V$, where $V(z)=\int_\Omega G(z,\zeta)\,d\mu_1(\zeta)$.
The Laplacian of $\psi$ is $\Delta\phi+\widetilde\mu_1-\mu_1=\widetilde\mu_1+\mu_2$ which has the desired 
properties. Moreover $|\phi-\psi|=|V|$ is bounded
by a constant times
\[
\int_{D(z,2R)}\left(1+\log \frac{2R}{d(z,\zeta)}\right) d\mu_1(\zeta).
\]

This integral is bounded by a constant times $\mu_1(D(z,2R))$, whenever $\mu_1$ 
a locally
doubling measure. This is Lemma~2.3 of \cite{Christ91}. 
The disk $D(z,2R)$ is covered by a 
bounded number of cubes $Q_j$, therefore the difference between 
$\psi$ and $\phi$
is bounded as claimed.
\end{proof}

\subsection{Proofs of Theorem~\ref{Christ} and \ref{meu}}
Let us start with Theorem~\ref{meu}. There are some weights  that are
particularly simple. These are the standard radial weights  $\phi(z)=\alpha\log
1/(1-|z|^2)$. The following lemma deals with this situation.

\begin{lema}\label{pes}
For any $\alpha\in (0,1)$ and $p\in[1,+\infty)$, the solution
\[
u(z)=\frac 1\pi \int_{\D}\frac{1-|\zeta|^2}{1-\bar\zeta
z}\frac{f(\zeta)}{z-\zeta}\, dm(\zeta)
\]
to the equation $\dbar u=f$ in $\D$ satisfies the estimate
\[
\int_\D |u(z)|^p (1-|z|)^{p\alpha-1} dm(z)\lesssim
\int_\D |f(z)(1-|z|)|^p (1-|z|)^{p\alpha-1} dm(z).
\]

Moreover,
\[
\sup_{\D} |u(z)|(1-|z|)^\alpha\lesssim \sup_{\D} |f(z)|(1-|z|)^{1+\alpha}.
\]
\end{lema}

\begin{proof}
This is an immediate consequence of H\"older's inequality, in the form of
Schur's test with respect to the measure $dm(z)/(1-|z|)$ for the kernel
$(1-|z|)^{\alpha}(1-|\zeta|)^{-\alpha}(1-|\zeta|^2)/(|1-\bar\zeta z||z-\zeta|)$.
(Replacing $\alpha$ by $\alpha/p$ one obtains the same estimate with the weight
$(1-|z|)^{\alpha-1}$ on both sides.)
\end{proof}

We take an arbitrary weight $\phi$ under the hypothesis of 
Theorem~\ref{meu}. By Lemma~\ref{regularitzacio} we may replace $\phi$ by an equivalent weight,
which only changes the constants in the estimates, so we may assume that
$(1-|z|^2)^2\Delta \phi>\varepsilon$ and $\Delta\phi$ is a locally
doubling measure with
respect to the hyperbolic measure; we may also assume that $\varepsilon\le 1$.
Consider the auxiliary subharmonic function
$\psi=\phi + (\varepsilon/8) \log (1-|z|^2)$. Since
$\Delta\log(1-|z|^2)=-4(1-|z|^2)^{-2}$, we have
$(1-|z|^2)^2\Delta \psi>\varepsilon/2$ 
and $\Delta\phi/2\le\Delta\psi\le\Delta\phi$, so that $\Delta\psi$ is still locally doubling. 
Using Theorem~\ref{multiplicador}, we can construct a
a holomorphic function $g$ such that  
$\frac{d(z,Z(g))^M}{\rho(z)^M}\lesssim |g|e^{-\psi}
\lesssim\frac{d(z,Z(g))}{\rho(z)}$.

To begin, let us assume that the support of $f$ is far from the zero set of the
multiplier $g$. That is, there is some $\delta>0$ such that 
$\frac{d(z,Z(g))}{\rho(z)}\ge\delta$.
Instead of solving the equation $\dbar u =f$, we consider the
auxiliary equation $\dbar v=f/g$. We take as a solution $v$ the function 
provided by Lemma~\ref{pes} (we take $\alpha=\varepsilon/8$). Then,
since $\dbar g=0$, the function $u=v g$ is a solution to $\dbar u =f$.
Moreover, because of Lemma~\ref{pes}, we know that for any $1\le p <\infty$
\[
\int_\D \frac{|u(z)/g(z)|^p}{(1-|z|)}(1-|z|)^{p\varepsilon/8} dm(z)\lesssim
\int_\D \frac{|f(z)/g(z)(1-|z|)|^p}{(1-|z|)}(1-|z|)^{p\varepsilon/8} dm(z).
\]

We always have that $|g|\lesssim e^{\psi}$, thus
\[
\int_\D \frac{|u(z)|^p}{(1-|z|)}e^{-p\phi(z)} dm(z)\lesssim
\int_\D \frac{|u(z)/g(z)|^p}{(1-|z|)}(1-|z|)^{p\varepsilon/8} dm(z),
\]
and since the support of $f$ is far from the zero sets of $g$, then
\[
\int_\D \frac{|f(z)/g(z)(1-|z|)|^p}{(1-|z|)}(1-|z|)^{p\varepsilon/8} dm(z)
\simeq
\int_\D \frac{|f(z)(1-|z|)|^p}{(1-|z|)}e^{-p\phi(z)} dm(z).
\]

The case $p=\infty$ follows with the same scheme.

Now, we must overcome the restriction on the support of $f$. 
We denote as above $\psi=\phi+\varepsilon/8 \log (1-|z|^2)$. For this
subharmonic function we take the set of
multipliers $h_i$ given by Proposition~\ref{tecnic} and its corresponding
exceptional sets $E_{h_i}$.

We split the domain into disjoint pieces: 
\[
\Omega=\bigl(\Omega\setminus E_{h_1}\bigr)\cup 
\bigl(E_{h_1}\setminus E_{h_2}\bigr)\cup
\bigl((E_{h_1}\cap E_{h_2})\setminus E_{h_3}\bigr)\cup\cdots\cup 
\bigl((E_{h_1}\cap\cdots\cap E_{h_{n-1}})\setminus E_{h_n}\bigr).
\] 
For the sake of simplicity we denote this partition 
 by $\Omega=\Omega_1\cup\cdots\cup\Omega_n$. In each $\Omega_i$ the
multiplier $|h_i|\simeq e^{\psi}$, and $|h_i|\lesssim e^\psi$ in all of
$\Omega$ (by the Remark after Theorem~\ref{multiplicador}). We can take as a solution to the equation
$\dbar u=f$ the function
\[
u(z)=\int_\D \left( \sum_{i=1}^n \frac{h_i(z)\chi_{\Omega_i(\zeta)}}{h_i(\zeta)}
\right)\frac 1\pi\frac{1-|\zeta|^2}{1-\bar\zeta
z}\frac{f(\zeta)}{z-\zeta}\, dm(\zeta)=\int_\D \kappa(z,\zeta)f(\zeta)
\, dm(\zeta).
\]
Thus, 
\[
|\kappa(z,\zeta)|\lesssim \frac{(1-|\zeta|^2)}{|1-\bar\zeta
z||\zeta-z|}\frac{(1-|z|)^{\varepsilon/8}e^{\phi(z)}}
{(1-|\zeta|)^{\varepsilon/8}e^{\phi(\zeta)}}.
\]
From this estimate and Schur's test, as in Lemma~\ref{pes} (with
$\alpha=\varepsilon/8$), the $L^p$ boundedness of the solution follows, for all
$p\in[1,\infty]$ and with the same solution $u$. This proves
Theorem~\ref{meu}.\qed

The same construction proves Theorem~\ref{Christ}. We  replace 
Lemma~\ref{pes} by the following one which is also a direct consequence of
H\"older's inequality:

\begin{lema}
For any $\alpha>0$ and $p\in[1,+\infty]$, the solution
\[
u(z)=\frac 1\pi \int_{\C}\frac{e^{2\alpha(\bar\zeta z -|\zeta|^2)}}{z-\zeta}
f(\zeta)\, dm(\zeta)
\]
to the equation $\dbar u=f$ in $\C$ satisfies the estimate 
$\|u(z)e^{-\alpha|z|^2}\|_p\lesssim \|f(\zeta)e^{-\alpha|\zeta|^2}\|_p$ for any
$p\in [1,\infty]$. 
\end{lema}

In this case, after replacing $\phi$ by an equivalent weight with
$\Delta\phi>\varepsilon\,dm$ (Lemma~\ref{regularitzacio}), the auxiliary subharmonic function $\psi$ is  $\phi-\varepsilon/8
|z|^2$, so that $\Delta\phi/2\le\Delta\psi\le\Delta\phi$. We take as a solution to the $\bar\partial$ equation the function
\[
\int_\C \left( \sum_{i=1}^n
\frac{h_i(z)\chi_{\Omega_i(\zeta)}}{h_i(\zeta)}\right)
 \frac 1\pi \frac{e^{(\varepsilon/4)(\bar\zeta z -|\zeta|^2)}}{z-\zeta}
f(\zeta)\,dm(\zeta)=\int_\C \kappa'(z,\zeta)f(\zeta)\, dm(\zeta).
\]
Therefore,
\[
|\kappa'(z,\zeta)|\lesssim  \frac{e^{\phi(z)}e^{-(\varepsilon/8)|z-\zeta|^2}}
{e^{\phi(\zeta)}|z-\zeta|}.
\]
This estimate proves Proposition~\ref{churrimangui} and Theorem~\ref{Christ}.
\qed

\subsection{The degenerate weight} 
We can prove this $\dbar$ estimate along the same lines . We need two
ingredients, a multiplier theorem and some $\dbar$ estimates when the weight
$\phi$ is of the form $\alpha |\Im z|$ for some $\alpha>0$. This is the
multiplier theorem that we need:

\begin{teor}\label{multreal}
Let $\phi$ be a subharmonic function in $\C$ such that the measure $\Delta
\phi$ is a locally 
doubling measure supported in the real line and  $\Delta \phi
(I(x,r))>1$ for some $r>0$ where $I(x,r)$ is the interval in $\R$ of center $x$
and radius $r$. There is a holomorphic function $f$ with zero set 
$\Lambda$
contained in $\R$ such that for any $\varepsilon>0$, $|f(z)|\simeq
e^{\phi(z)}$, for all $z$ such that 
$|z-\lambda_n|\ge \varepsilon \rho(\lambda_n)$ for all
$\lambda_n\in Z(f)$.
\end{teor}

\begin{proof} 
The proof of this theorem is the same as in Theorem~\ref{multiplicador} when
$\Omega=\C$,  except that at some points it is easier. For instance, it is
trivial to split the real line into intervals $I_k=[x_k,x_{k+1})$ all of mass
$4\pi$, because $x\mapsto\Delta\phi([0,x))$ is continuous and strictly
increasing; these intervals have length $|I_k|\simeq\rho(x)$ for $x\in I_k$.
Moreover, on the line it is enough to match the moments of order $0$ and $1$,
so Lemma~\ref{moments} is not needed and the zeros can be taken real (the
points given by Lemma~\ref{moments} need not be real): we put a double
zero at the barycenter $\lambda_k=(4\pi)^{-1}\int_{I_k}t\,d\Delta\phi(t)$ of each
interval, which lies in $I_k$ at distance $\gtrsim|I_k|$ from its endpoints by
the doubling property. If $|z-\lambda_k|\ge 2|I_k|$, Taylor's formula of order
$2$ gives
\[
\left|\int_{I_k}\log|z-t|\,d\Delta\phi(t)-4\pi\log|z-\lambda_k|\right|\lesssim
\frac{|I_k|^2}{|z-\lambda_k|^2},
\]
and the sum of these errors is bounded uniformly in $z$: the intervals at
distance $\simeq d$ from $z$ have total length $\lesssim d$ and, by
Lemma~\ref{dis-mes} (applied twice), $|I_k|/d\lesssim(\rho(\Re z)/d)^{\gamma^2}$ when $d\le 1$,
while $|I_k|\lesssim 1$. The intervals close to $z$ are treated as in
Theorem~\ref{multiplicador}.
\end{proof}

On the other hand the $\dbar$-estimate that we need in the flat case, i.e. when
$\phi=\alpha |\Im z|$, is not as easy as in the disk or the plane; we need the
following theorem, a proof of which can be found in \cite{OrtSei99}:

\begin{teor*}
Consider the equation $\dbar u=\mu$, where $\mu$ is a compactly supported
measure such that $e^{-\alpha|\Im z|}d|\mu|$ is
a two-sided  Carleson measure for some $\alpha>0$. 
Then there is a solution $u$ with
\[
\limsup_{z\to\infty} |u(z)|e^{-\alpha|\Im z|}=0 \quad \text{ and } \quad
|u(x)|\le C \left(1+ \int_{|z-x|<1}\frac {d|\mu|(z)}{|x-z|}\right)
\]
for any $x\in \R$, where $C$ only depends on $\alpha$ and on the Carleson constant of 
$e^{-\alpha|\Im z|}d|\mu|$. 
\end{teor*}
(In \cite[Lemma~4]{OrtSei99} the first condition is stated as
$\limsup_{z\to\infty}|u(z)|e^{-\alpha|\Im z|}<\infty$, but the solution
constructed there satisfies it with $=0$, which is what we need.)

These two ingredients together prove Theorem~\ref{degenerat} in the same way 
as we proved Theorem~\ref{Christ} and Theorem~\ref{meu}, with two
modifications. First, we need multipliers for the weight $\phi-\alpha|\Im z|$
for some $\alpha>0$, whose Laplacian $\Delta\phi-2\alpha\,dx$ is a signed
measure on $\R$. The proof of Theorem~\ref{multreal} applies to it without
changes when $0<\alpha\le\alpha_0$ is small enough: the points $x_k$ are
chosen with $(\Delta\phi-2\alpha\,dx)([x_k,x_{k+1}))=4\pi$, and the intervals
$I_k$ still have length $\simeq\rho$ and barycenters (with respect to
$\Delta\phi-2\alpha\,dx$) at distance $\gtrsim|I_k|$ from their endpoints.
Second, replacing each double zero $\lambda_k$ by the two simple zeros
$\lambda_k\pm s|I_k|$, with $s>0$ small, we obtain a second multiplier whose
exceptional set is disjoint from that of the first one (as in
Proposition~\ref{tecnic}). In this way we get two entire functions $f_1,f_2$
with real zeros and a partition $\C=\Omega_1\cup\Omega_2$ into Borel sets such
that $|f_i(z)|\lesssim e^{\phi(z)-\alpha|\Im z|}$ in $\C$ and
$|f_i(z)|\simeq e^{\phi(z)-\alpha|\Im z|}$ in $\Omega_i$. Let $v_i$ be the
solution of $\dbar v_i=\chi_{\Omega_i}\mu/f_i$ given by the theorem above,
which applies because $e^{-\alpha|\Im z|}\chi_{\Omega_i}d|\mu|/|f_i|\lesssim
e^{-\phi}d|\mu|$. Then $u=f_1v_1+f_2v_2$ solves $\dbar u=\mu$,
$|u(z)|e^{-\phi(z)}\lesssim\sum_i|v_i(z)|e^{-\alpha|\Im z|}\to 0$ as
$z\to\infty$ and, since $|f_i(x)|e^{-\phi(x)}\lesssim 1$ for $x\in\R$ and
$1/|f_i(z)|\lesssim e^{\alpha}e^{-\phi(z)}$ for $z\in\Omega_i$ with
$|\Im z|<1$,
\[
|u(x)|e^{-\phi(x)}\lesssim\sum_{i=1}^2|v_i(x)|\lesssim
1+\int_{|z-x|<1}\frac{e^{-\phi(z)}\,d|\mu|(z)}{|x-z|}.
\]
This is why the weight $e^{-\phi(z)}$ appears inside the integral in
Theorem~\ref{degenerat}.
\subsection*{Added in proof} Richard Rochberg has informed me of un 
unpublished manuscript of
Tom Wolff from 1988 where some of the ideas concerning the multiplier
(Theorem~\ref{multiplicador}) are already present.

\end{document}